\documentstyle[11pt, amssymb]{article}

\begin{document}

\begin{flushleft}
{\Large {\bf Higher genus polylogarithms on families of Riemann surfaces}}
\end{flushleft}

\begin{flushleft}
{\large {\bf Takashi Ichikawa}} 
\end{flushleft}

\begin{flushleft} 
Department of Mathematics, Faculty of Science and Engineering, 
Saga University, Saga 840-8502, Japan. E-mail:  ichikawn@cc.saga-u.ac.jp 
\end{flushleft} 

%\begin{flushleft}
%{\it MSC:} 11G55, 11M32, 14H10, 14H15, 14H30, 20F34, 57K20  
%\end{flushleft}

%\begin{flushleft}
%{\it Keywords:} Moduli space of algebraic curves, Universal Mumford curve, 
%Elliptic associator, Polylogarithm, Multiple zeta value, Teichm\"{u}ller groupoid 
%\end{flushleft}

\noindent
{\bf Abstract:} 
We construct polylogarithms on families of pointed Riemann surfaces of any genus 
which describe monodromies of meromorphic connections with simple poles. 
Furthermore, we show that the polylogaritms are computable as power series 
in deformation parameters and their logarithms associated with the families. 
%\vspace{4ex}

%\noindent
%MSC-Class: 11G20, 11G55, 11M32, 14H10, 14H15, 14D15, 32G20 
%\vspace{4ex}

%\noindent
%\hrulefill
%\vspace{2ex}

\section{Introduction} 

{\it Polylogarithm functions}, or {\it polylogarithms} for short, 
describe unipotent periods which appear as monodromies of nilpotent connections 
defined on families of algebraic varieties. 
Especially, polylogarithms on families of pointed Riemann surfaces of genus $0$, 
called rational polylogarithms, 
become important objects for studying Feynman integrals in quantum field theory and 
motives in arithmetic geometry. 
A good reference on this subject is \cite{GiEG}. 
Furthermore, polylogarithms on Riemann surfaces of positive genus and 
associated Feynman integrals were studied 
(see \cite{Brod, Bo, W, Bou} and references there in). 
Very recently, 
higher genus polylogaritms were presented by D'Hoker-Hidding-Schlotterer \cite{DHS} 
which are considered as monodromies of a single-valued, 
but non-meromorphic connection on each Riemann surface of any genus. 

In this paper, we construct higher genus polylogaritms on families of 
pointed Riemann surfaces which describe monodromies of meromorphic connections 
with simple poles (but not single-valued), 
and show that the polylogaritms are computable as power series 
in deformation parameters and their logarithms associated with the families. 
First, we construct a polylogarithm sheaf on a special family of Riemann surfaces 
by gluing the KZ (Knizhnik-Zamolodchikov) connections and its elliptic extensions 
known as the elliptic KZB (Knizhnik-Zamolodchikov-Bernard) connections 
by Calaque-Enriquez-Etingof \cite{CEE} (see also \cite{E1, H}). 
Second, we show that this polylogarithm sheaf can be analytically continued to 
general families of pointed Riemann surfaces, 
and that the associated polylogarithms are described by the power series as above 
whose coefficients are essentially expressed by mutiple zeta values. 
Our result is an extension of results of Deligne \cite{D2}, Brown \cite{Brow} and 
Banks-Panzer-Pym \cite{BPP} in the genus $0$ case, 
and is useful even in computing polylogarithms on elliptic curves with multiple points 
which was considered in \cite{Bo}.

\section{Polylogarithms and associators}

\subsection{Polylogarithms for Riemann surfaces} 

Fix non-negative integers $g, n$ such that $2g - 2 + n > 0$, 
let $({\mathcal R}; p_{1},..., p_{n})$ be a family over a complex space $S$ of 
$n$-pointed compact Riemann surfaces of genus $g$, 
and put ${\mathcal R}^{\circ} = {\mathcal R} - \{ p_{1},..., p_{n} \}$. 
Then a {\it polylogarithm sheaf} on $({\mathcal R}; p_{1},..., p_{n})$ consists of 
a vector bundle ${\mathcal V}$ on ${\mathcal R}$ with flat connection ${\mathcal F}$ 
such that ${\mathcal F}$ is holomorphic on ${\mathcal R}^{\circ}$ 
and it has simple poles along $p_{i}$ $(i = 1,..., n)$. 
Then the associated {\it unipoent periods} are defined as the monodromies 
between (tangential) sections of ${\mathcal R}^{\circ}/S$ 
which are described by multi-valued functions on $S$ 
called {\it polylogarithms} on $({\mathcal R}; p_{1},..., p_{n})/S$. 
We consider the case that ${\mathcal R}$ is the projective line, 
namely the Riemann sphere ${\mathbb P}^{1} = {\mathbb C} \cup \{ \infty \}$, 
or an elliptic curve.

\subsection{KZ connection and associator} 

For an integer $n \geq 3$,  
a {\it KZ connection} (Knizhnik-Zamolodchikov) on an $n$-pointed projective line 
$({\mathbb P}^{1}; p_{1},..., p_{n})$ is defined as a trivial bundle with flat connection 
which has simple poles at $p_{i}$ with residue $X_{p_{i}}$ $(i = 1,..., n)$, 
where $X_{p_{i}}$ are symbols satisfying $\sum_{i=1}^{n} X_{p_{i}} = 0$. 
For KZ connections on pointed projective lines 
$P = ({\mathbb P}^{1}; p_{1},..., p_{n})$ and $Q = ({\mathbb P}^{1}; q_{1},..., q_{m})$ 
such that $X_{p_{n}} + X_{q_{m}} = 0$, 
one obtain an $(n + m - 2)$-pointed projective line $P \cup Q$ by gluing $P, Q$ 
around $p_{n}, q_{m}$ via the relation that $(z - p_{n})(w - q_{m})$ is a non-zero constant, 
where $z, w$ are fixed coordinates on $P, Q$ respectively 
satisfying $p_{n}, q_{m} \neq \infty$ for sake of simplicity. 
Then a connection on $P \cup Q$ obtained by gluing the KZ connections 
on $P, Q$ via the relation $X_{p_{n}} dz/(z - p_{n}) = X_{q_{m}} dw/(w - q_{m})$ 
becomes a KZ connection since 
$\sum_{i = 1}^{n-1} X_{p_{i}} + \sum_{j = 1}^{m-1} X_{q_{j}} = 0$. 

We take a coordinate $z$ on ${\mathbb P}^{1}$ such that $p_{n} = \infty$. 
Then the connection form of a KZ connection on $P$ is given by the equation 
$$
d - \sum_{i=1}^{n-1} X_{p_{i}} \frac{dz}{z - p_{i}} = 0, 
$$
and hence the associated monodromy along a path 
$\gamma \subset {\mathbb P}^{1}$ from $a$ to $b$ is 
a sum of words of $X_{p_{1}},... X_{p_{n}}$ whose coefficients are 
{\it rational polylogarithms} given by iterated integrals 
$$
\int_{\gamma} \frac{dz}{z - p_{i(1)}} \cdots \frac{dz}{z - p_{i(k)}}, 
$$
where $p_{i(j)} \in \{ p_{1},..., p_{n-1} \}$ (cf. \cite{Brow, GoM}). 
These polylogarithms are convergent when $a \neq p_{i(1)}$, $b \neq p_{i(k)}$, 
and are generally described in the shuffle regularization (cf. \cite[2.1.5]{BPP}) 
as sums of convergent polylogarithms times 
$$
\frac{1}{j!} \left( \int_{\gamma} \frac{dz}{z - a} \right)^{j}   
\frac{1}{k!} \left( \int_{\gamma} \frac{dz}{z - b} \right)^{k} 
$$ 
when $a = p_{i(1)}$ or $b = p_{i(k)}$. 

We assume that $n = 3$, 
and hence one can take $p_{1}, p_{2}$ as $0, 1$ respectively. 
Then tangential points on $({\mathbb P}^{1}; 0, 1, \infty)$ are defined as six 
unit and real tangent vectors based at either $0, 1$ or $\infty$. 
For example, 
$\overrightarrow{01}$ (resp. $\overrightarrow{10}$) denote the tangential points 
at $0$ (resp. $1$) with positive (resp. negative) orientations. 
Then there exists a unique solution $G_{0}(z)$ of the KZ connection 
$$
d - X_{0} \frac{dz}{z} - X_{1} \frac{dz}{z - 1} = 0 
\eqno(2.1)
$$
which is normalized at $\overrightarrow{01}$, 
namely satisfying that 
$$
\lim_{z \downarrow 0} \frac{G_{0}(z)}{z^{X_{0}}} = 1; \ \ 
z^{X_{0}} := \sum_{k=0}^{\infty} \frac{\log(z)^{k}}{k!} (X_{0})^{k} \ (0 < z < 1). 
$$
Denote by ${\mathbb C} \left\langle\!\left\langle X_{0}, X_{1} \right\rangle\!\right\rangle$ 
the ring of non-commutative formal power series over ${\mathbb C}$ 
in the symbols $X_{0}$ and $X_{1}$. 
Then $G_{0}(z) \cdot z^{-X_{0}}$ is an element of 
${\mathbb C} \left\langle\!\left\langle X_{0}, X_{1} \right\rangle\!\right\rangle$ 
whose constant term is $1$ and coefficients are rational polylogarithm functions 
\begin{eqnarray*}
{\rm Li}_{k_{1},..., k_{l}}(z) 
& = & 
\int_{0}^{z} \underbrace{w_{0} \cdots w_{0}}_{k_{l}-1} 
w_{1} \underbrace{w_{0} \cdots w_{0}}_{k_{l-1}-1} w_{1} \cdots w_{1} 
\underbrace{w_{0} \cdots w_{0}}_{k_{1}-1} w_{1} 
\\
& = & 
\sum_{0 < n_{1} < \cdots < n_{l}} \frac{z^{n_{l}}}{n_{1}^{k_{1}} \cdots n_{l}^{k_{l}}} 
\in {\mathbb Q}[[z]] \ \ (|z| < 1), 
\end{eqnarray*}
where $w_{0} = dz/z$, $w_{1} = dz/(1-z)$. 
Furthermore, 
denote by $G_{1}(z)$ the solution of (2.1) normalized at $\overrightarrow{10}$, i.e.,  
$$
\lim_{z \uparrow 1} \frac{G_{1}(z)}{(1 - z)^{X_{1}}} = 1. 
$$
Then the unipotent period of (2.1) from $\overrightarrow{01}$ to $\overrightarrow{10}$  becomes the connection matrix 
$$
\Phi_{\rm KZ}(X_{0}, X_{1}) = G_{1}(z)^{-1} \cdot G_{0}(z) 
$$ 
which is given by the {\it KZ associator} as an element of 
${\mathbb C} \left\langle\!\left\langle X_{0}, X_{1} \right\rangle\!\right\rangle$ 
with coefficients expressed by {\it multiple zeta values} 
$$
\zeta({k_{1},..., k_{l}}) = 
\sum_{0 < n_{1} < \cdots < n_{l}} \frac{1}{n_{1}^{k_{1}} \cdots n_{l}^{k_{l}}} 
$$
for $k_{l} > 1$.

\subsection{Reduced elliptic KZB connection} 

The universal elliptic KZB (Knizhnik-Zamolodchikov-Bernard) connection 
was constructed in \cite[Theorem 12]{CEE} as a vector bundle with flat connection 
over a family of pointed elliptic curves. 
Following \cite{H}, 
we review its reduction to the elliptic curves 
$$
E_{q} = \left( {\mathbb P}^{1} - \{ 0, \infty \} \right) / \langle q \rangle; \ 
\langle q \rangle = \left\{ q^{a} \ | \ a \in {\mathbb Z} \right\} 
$$
with marked point $1$ for $0 < |q| < 1$ which are expressed as 
${\mathbb C}/({\mathbb Z} + {\mathbb Z} \tau)$ if $q =  e^{2 \pi i \tau}$. 
Let 
$$
\theta (u, \tau) = \sum_{n \in {\mathbb Z}} (-1)^{n} 
q^{\frac{1}{2} \left( n + \frac{1}{2} \right)^{2}} e^{\left( n + \frac{1}{2} \right) u}; \
q = e^{2 \pi i \tau} 
$$ 
be the classical theta function, 
put 
$$
F \left( \xi, \eta, \tau \right) = 
2 \pi i \frac{ \, \theta' (0, \tau) \, \theta (2 \pi i (\xi + \eta), \tau) \, }
{\theta (2 \pi i \xi, \tau) \, \theta (2 \pi i \eta, \tau)} 
$$
and denote by 
${\mathbb C} \left\langle\!\left\langle {\bf t}, {\bf a} \right\rangle\!\right\rangle$ 
the ring of non-commutative formal power series over ${\mathbb C}$ 
in the symbols ${\bf t}$ and ${\bf a}$. 
Then putting $d \tau = 0$ on the universal elliptic KZB equation \cite[9.2]{H}, 
we have the reduced elliptic KZB equation on 
$E_{q} = {\mathbb C}/({\mathbb Z} + {\mathbb Z} \tau)$ with fibers 
${\mathbb C} \left\langle\!\left\langle {\bf t}, {\bf a} \right\rangle\!\right\rangle$ 
which is defined as 
$$
d + {\bf t} F \left( \xi, {\bf t}, \tau \right) ({\bf a}) \ d \xi = 0, 
\eqno(2.2) 
$$
where $f({\bf t})({\bf a}) = f \left( {\rm ad}_{\bf t} \right)({\bf a})$. 
Put $T = 2 \pi i {\bf t}, A = (2 \pi i)^{-1} {\bf a}$ and $q = e^{2 \pi i \tau}, z = e^{2 \pi i \xi}$. 
By \cite[(8.1)]{H}, 
\begin{eqnarray*} 
F \left( \xi, {\bf t}, \tau \right) 
& = & 
\pi i \left[ \coth(\pi i \xi) + \coth(\pi i {\bf t}) \right] + 
4 \pi \sum_{n=1}^{\infty} 
\left(\sum_{d|n} \sin \left[ 2 \pi \left( \frac{n}{d} \xi + d {\bf t} \right) \right] \right) q^{n} 
\\ 
& = & 
\pi i \left( \frac{z + 1}{z - 1} + \frac{e^{T} + 1}{e^{T} - 1} \right) - 
2 \pi i \sum_{n=1}^{\infty} \sum_{d|n} 
\left( z^{n/d} e^{dT} - z^{-n/d} e^{-dT} \right) q^{n}, 
\end{eqnarray*} 
and hence 
\begin{eqnarray*} 
& & 
{\bf t} F \left( \xi, {\bf t}, \tau \right)({\bf a}) d \xi 
\\ 
& = & 
[T, A] \frac{dz}{z - 1} 
+ \left( \frac{T}{e^{T} - 1} \right)(A) \frac{dz}{z} 
- T \sum_{n=1}^{\infty} \sum_{d|n} 
\left( z^{n/d} e^{dT} - z^{-n/d} e^{-dT} \right)(A) \, q^{n} \frac{dz}{z}. 
\end{eqnarray*} 
Therefore, 
the reduced elliptic KZB connection (2.2) over $E_{q}$ with fibers 
${\mathbb C} \left\langle\!\left\langle T, A \right\rangle\!\right\rangle$ 
has a simple pole at $z = 1$ with residue $W_{1} := - [T, A]$, 
and has poles at $z = 0, \infty$ with residues 
$$
W_{0} := - \left( \frac{T}{e^{T} - 1} \right) (A), \ 
W_{\infty} := -(W_{0} + W_{1}) = \left( \frac{-T}{e^{-T} - 1} \right) (A) = e^{T}(- W_{0}) 
$$
respectively. 
Furthermore, 
the restriction of (2.2) to $E_{q}|_{q = 0}$ becomes the KZ equation (2.1), 
where $X_{i} = W_{i}$ $(i = 1, 2)$ (cf. \cite[Section 12]{H}).

\subsection{Elliptic multiple zeta values}

Tangential points on $E_{q}$ are defined as two unit and real tangent vectors based at $1$. 
Then the connection matrices of (2.2) between these tangential points are expressed by 
$A(\tau)$, $B(\tau)$ which are defined in \cite[6.2]{E1} and studied in \cite{E2}.  
Especially, it is shown in \cite[Section 5]{E2} that 
the associated polylogarithms are expressed as {\it elliptic multiple zeta values} 
$$
\sum_{n=0}^{\infty} I_{n} q^{n}, \ \ 
\sum_{m \in {\mathbb Z}} \sum_{n=0}^{\infty} J_{m, n} \tau^{m} q^{n}, 
$$
where $I_{n}$ and $J_{m, n}$ are computable ${\mathbb Q}$-linear sums of 
products of $\pi i$ and multiple zeta values.

\section{Higher genus polylogarithm}

\subsection{Polylogarithm sheaf} 

We construct polylogarithm sheaves. 
Let $\Delta_{0} = (V_{0}, E_{0}, T_{0})$ be a stable graph with orientation of $(g, n)$-type 
which is {\it trivalent} in the sense that each vertex has just three branches. 
Furthermore, we assume that $E_{0}$ contains $g$ loops denoted by $e_{1},..., e_{g}$, 
namely $\Delta_{0}$ consists of its unique maximal subtree $M_{0}$ and $e_{i}$ 
connected by vertices $v_{e_{i}} = v_{-e_{i}}$ $(i = 1,..., g)$. 
Consider the ring 
$$
{\mathcal X}_{g, n} = 
{\mathbb C} \left\langle\!\left\langle X_{h}, T_{i}, A_{i} \right\rangle\!\right\rangle 
$$
of non-commutative formal power series over ${\mathbb C}$ in the symbols 
$X_{h}$ $(h \in \pm M_{0} \cup T_{0})$ and $T_{i}, A_{i}$ $(i = 1,..., g)$ 
satisfying the conditions 
\begin{itemize}

\item 
For each $e \in E_{0} - \{ e_{1},..., e_{g} \}$, $X_{e} + X_{-e} = 0$. 

\item 
For each $v \in V$, 
the sum of $X_{h}$ for $h \in \pm E_{0} \cup T_{0}$ with $v_{h} = v$ is equal to $0$, 
where 
\begin{eqnarray*}
X_{e_{i}} 
& = & 
\left( \frac{- T_{i}}{e^{T_{i}} - 1} \right) (A_{i}), 
\\ 
X_{-e_{i}} 
& = & 
\left( \frac{- T_{i}}{e^{-T_{i}} - 1} \right) (A_{i}) = - e^{T_{i}} \left( X_{e_{i}} \right). 
\end{eqnarray*} 

\end{itemize} 
Since $\Delta_{0}$ is trivalent, as in A.2, 
let ${\mathcal R}_{\Delta_{0}}$ be the deformation of the union of 
$$
P_{v} = ({\mathbb P}^{1}; 0, 1, \infty) \ \ (v \in V_{0})
$$ 
by the relations $\xi_{h} \cdot \xi_{-h} = y_{h} \ \ (h \in \pm E_{0})$, 
where $\xi_{h}$ denotes the local coordinate at $x_{h} \in P_{v_{h}}$. 
Denote by $\left( {\mathcal V}_{v}, {\mathcal F}_{v} \right)$ the KZ connections on $P_{v}$ 
$(v \in V_{0} - \{ v_{e_{1}},..., v_{e_{g}} \})$,  
and by $\left( {\mathcal V}_{i}, {\mathcal F}_{i} \right)$ 
the reduced elliptic KZ connections on 
$$ 
E_{i} = \left( ({\mathbb P}^{1} - \{ 0, 1 \})/\langle y_{e_{i}} \rangle; 1 \right) 
\ \ (i = 1,..., g)
$$ 
which have simple poles at $1$ with residue $- [T_{i}, A_{i}]$. 
Then by the above conditions, 
one can glue $\left( {\mathcal V}_{v}, {\mathcal F}_{v} \right)$ $(v \neq v_{e_{1}},..., v_{e_{g}})$  
and $\left( {\mathcal V}_{i}, {\mathcal F}_{i} \right)$ $(i = 1,..., g)$ 
around $y_{e} = 0$ $(e \neq e_{1},..., e_{g})$. 
Therefore, 
we obtain a polylogarithm sheaf 
$\left( {\mathcal V}_{\Delta_{0}}, {\mathcal F}_{\Delta_{0}} \right)$ 
on ${\mathcal R}_{\Delta_{0}}$ as a vector bundle with flat connection 
which is meromorphic with simple poles at $x_{t}$ $(t \in T)$. 

Let $\Delta = (V, E, T)$ be a stable graph of $(g, n)$-type with orientation. 
Then $\Delta$ is obtained from $\Delta_{0}$ by repeating alterations given in A.2 
{\it not contracting $e_{1},..., e_{g}$ to points.} 
We attach symbols $X^{\Delta}_{h}$ with each $h \in \pm E \cup T$ 
by the following rules: 
\begin{itemize}

\item 
If $\Delta = \Delta_{0}$, 
then $X^{\Delta}_{h} = X_{h}$. 

\item 
For any $v \in V$, 
the sum of $X^{\Delta}_{h}$ $(h \in \pm E \cup T$ with $v_{h} = v)$ is $0$. 

\item 
Let $\Delta = (V, E, T)$ and $\Delta' = (V', E', T')$ be stable graphs of $(g, n)$-type 
which are related as in A.2. 
Then $X^{\Delta}_{h} = X^{\Delta'}_{h}$ for $h \neq \pm h_{0}$. 

\end{itemize}
Then $({\mathcal V}_{\Delta_{0}}, {\mathcal F}_{\Delta_{0}})$ on ${\mathcal R}_{\Delta_{0}}$ 
can be analytically continued to a vector bundle with flat connection 
$({\mathcal V}, {\mathcal F})$ on a family of $n$-pointed Riemann surfaces of genus $g$ 
close to degenerate complex curves which we call the {\it universal polylogarithm sheaf}.

\subsection{Higher genus polylogarithms} 

For a stable graph $\Delta$ of $(g, n)$-type, 
denote by $\left( {\mathcal R}_{\Delta}; x_{t} \ (t \in T) \right)$ 
the associated family of $n$-pointed Riemann surfaces of genus $g$ as in A.2 
under identifying $T$ with $\{ 1,..., n \}$. 
Then {\it (multiple) polylogarithms} on ${\mathcal R}_{\Delta}$ are defined as 
the coefficients of words of $X^{\Delta}_{h}$ $(h \in \pm E \cup T)$ which give rise to 
the unipotent periods between tangential points at $x_{t} \in P_{v_{t}}$ 
of the universal polylogarithm sheaf. 

First, we compute polylogarithms on ${\mathcal R}_{\Delta_{0}}$, 
where $\Delta_{0} = (V_{0}, E_{0}, T_{0})$ is the trivalent graph of $(g, n)$-type given in 3.1. 
We consider the unipotent periods of 
$\left( {\mathcal V}_{\Delta_{0}}, {\mathcal F}_{\Delta_{0}} \right)$ 
obtained as products of its monodromies along 
\begin{itemize}

\item[(1)] 
paths between tangential points on $P_{v}$ for each $v \in V_{0} - \{ v_{e_{1}},..., v_{e_{g}} \}$. 

\item[(2)] 
paths between tangential points at $x_{e}$ and $x_{-e}$ 
for each $e \in E_{0} - \{ e_{1},..., e_{g} \}$. 

\item[(3)] 
closed paths in $E_{i}$ $(i = 1,..., g)$ between tangential points at $1$. 

\end{itemize}
The monodromies along paths in (1) are products of 
$$
e^{\pi i X^{\Delta_{0}}_{h}} = 
\sum_{k=0}^{\infty} \frac{(\pi i)^{k}}{k !} \left( X^{\Delta_{0}}_{h} \right)^{k} 
\ \ (v_{h} = v), 
$$ 
and the KZ associators 
$$
\Phi_{\rm KZ} \left( X^{\Delta_{0}}_{h}, X^{\Delta_{0}}_{h'} \right) 
\ \ (v_{h} = v_{h'} = v) 
$$
which are seen in 2.2 to be non-commutative formal power series 
in $X_{h}^{\Delta_{0}}, X_{h'}^{\Delta_{0}}$ 
whose coefficients are expressed by $\pi i$ and multiple zeta values. 
The monodromies along paths in (2) are 
$$
\xi_{h}^{X^{\Delta_{0}}_{h}} \cdot \left( \xi_{-h}^{X^{\Delta_{0}}_{-h}} \right)^{-1} = 
y_{h}^{X^{\Delta_{0}}_{h}} = 
\sum_{k=0}^{\infty} \frac{\log (y_{h})^{k}}{k!} \left( X^{\Delta_{0}}_{h} \right)^{k} 
\ \ \mbox{around $y_{h} = 0$}. 
$$
Furthermore, as is seen in 2.4, the associated polylogarithms for paths in (3) are 
$$
\sum_{n=0}^{\infty} I_{n} y_{e_{0}}^{n}, \ \ 
\sum_{m \in {\mathbb Z}} \sum_{n=0}^{\infty} 
\frac{J_{m, n}}{(2 \pi i)^{m}} \log(y_{e_{0}})^{m} y_{e_{0}}^{n}, 
$$
where $I_{n}, J_{m, n}$ are computable ${\mathbb Q}$-linear sums of 
products of $\pi i$ and multiple zeta values. 
Therefore, we can compute the products of the above monodromies, 
and the associated polylogarithms. 

Second, 
we study the analytic continuation of the monodromies 
$$
\xi_{h}^{X^{\Delta}_{h}} \cdot \left( \xi_{-h}^{X^{\Delta}_{-h}} \right)^{-1} = y_{h}^{X^{\Delta}_{h}} 
$$
along paths in (2) under a fusing (associative) move 
which corresponds to the move of $y_{h}$ from $\overrightarrow{01}$ 
to $\overrightarrow{10}$. 
Take two $3$-pointed projective lines $P_{i} = ({\mathbb P}^{1}; 0, 1, \infty)$ $(i = 1, 2)$, 
and glue them around $\infty \in P_{1}$ and $0 \in P_{2}$ by the relation (A.1), 
i.e., $z_{2} = \phi(z_{1}) = y z_{1}$, 
where $z_{i}$ are the coordinates on $P_{i}$ and $y = y_{h}$ denotes the deformation parameter. 
Since $\phi(0) = 0$ and $\phi(1) = y$, 
we have a projective line with $4$ marked points $0, 1, \infty$ and $y$, 
and hence for the symbol $X = X^{\Delta}_{h}$, 
$$
y^{X} = \sum_{k=0}^{\infty} \frac{\log(y)^{k}}{k!} X^{k} \ \ 
\mbox{around $y = 0$}
$$
is analytically continued to 
$$
\left( 1 - (1 - y) \right)^{X} = 
\sum_{k=0}^{\infty} \frac{X (X-1) \cdots (X - k + 1)}{k!} (-1)^{k} (1 - y)^{k} \ \ 
\mbox{around $y = 1$}. 
$$ 
Since all trivalent graphs of $(g, n)$-type are connected by fusing moves, 
if a stable graph $\Delta$ of $(g, n)$-type is trivalent, 
then there exist explicit formulas of polylogarithms on ${\mathcal R}_{\Delta}$. 
Furthermore, by the consideration in 2.2, 
we have a similar result to \cite[Theorem 2.21]{BPP} (cf. \cite{Brow}) 
for the polylogarithms when $\Delta$ is not trivalent.

\appendix

\section{Families of Schottky uniformized Riemann surfaces}  

\subsection{Schottky uniformized deformation} 

A Schottky group $\Gamma$ of rank $g$ is defined as a free group with generators 
$\gamma_{i}$ $(i = 1,..., g)$ in 
$PGL_{2}({\mathbb C}) = GL_{2}({\mathbb C})/({\mathbb C} - \{ 0 \})$ 
which map Jordan curves  
$C_{i} \subset {\mathbb P}^{1} = {\mathbb C} \cup \{ \infty \}$ 
to other Jordan curves $C_{-i} \subset {\mathbb P}^{1}$ 
with orientation reversed, 
where $C_{\pm 1},..., C_{\pm g}$ with their interiors are mutually disjoint. 
Each element $\gamma \in \Gamma - \{ 1 \}$ is conjugate 
to an element of $PGL_{2}({\mathbb C})$ sending $z$ to $\beta_{\gamma} z$ 
for some $\beta_{\gamma} \in {\mathbb C}$ with $0 < |\beta_{\gamma}| < 1$ 
which is called the {\it multiplier} of $\gamma$. 
Therefore, we have 
$$
\frac{\gamma(z) - \alpha_{\gamma}}{z - \alpha_{\gamma}} = 
\beta_{\gamma} \frac{\gamma(z) - \alpha'_{\gamma}}{z - \alpha'_{\gamma}} 
$$
for some element $\alpha_{\gamma}, \alpha'_{\gamma}$ of ${\mathbb P}^{1}$ 
called the {\it attractive}, {\it repulsive} fixed points of $\gamma$ 
respectively. 
Then the discontinuity set $\Omega_{\Gamma} \subset {\mathbb P}^{1}$ 
under the action of $\Gamma$ has a fundamental domain given by 
the complement of the union of the interiors of $C_{\pm 1},..., C_{\pm g}$. 
The quotient space $R_{\Gamma} = \Omega_{\Gamma}/\Gamma$ 
is a (compact) Riemann surface of genus $g$ 
which is called {\it Schottky uniformized} by $\Gamma$. 
Furthermore, by a result of Koebe, 
every Riemann surface of genus $g$ can be Schottky uniformized, 
and an isomorphism between two Riemann surfaces gives an isomorphism 
between the corresponding Schottky groups given by the conjugation of an element. 

We construct families of Schottky uniformized pointed Riemann surfaces 
as universal deformations of degenerate pointed complex curves which are stable curves 
(cf. \cite{DM, K, KM}) such that the normalizations of their irreducible components 
are projective lines. 
Let $\Delta = (V, E, T)$ be a connected graph which is a collection of 
three finite sets $V$ of vertices, $E$ of edges, $T$ of tails and two boundary maps 
$$
b : T \rightarrow V, 
\ \ b : E \longrightarrow \left( V \cup \{ \mbox{unordered pairs of elements of $V$} \} \right). 
$$
We assume that a graph $\Delta = (V, E, T)$ is {\it stable}, 
namely its each vertex has at least three branches, 
and call it {\it of $(g, n)$-type} if $\dim H_{1}(\Delta) = g$ 
and the number of elements of $T$ is $n$. 
An {\it orientation} of a stable graph $\Delta = (V, E, T)$ means 
giving an orientation of each $e \in E$. 
Under an orientation of $\Delta$, 
denote by $\pm E = \{ e, -e \ | \ e \in E \}$ the set of oriented edges, 
and by $v_{h}$ the terminal vertex of $h \in \pm E$ (resp. the boundary vertex of $h \in T)$. 
For each $h \in \pm E$, 
denote by let $| h | \in E$ be the edge $h$ without orientation. 

For a stable graph $\Delta = (V, E, T)$ of $(g, n)$-type with orientation, 
take a subset ${\mathcal E}$ of $\pm E \cup T$ 
whose complement ${\mathcal E}_{\infty}$ satisfies the condition that 
$$
\pm E \cap {\mathcal E}_{\infty} \cap \{ -h \ | \ h \in {\mathcal E}_{\infty} \} 
\ = \ 
\emptyset, 
$$ 
and that $v_{h} \neq v_{h'}$ for any distinct $h, h' \in {\mathcal E}_{\infty}$. 
We attach Riemann spheres $P_{v} = {\mathbb P}^{1}$ for $v \in V$, 
and complex moduli parameters $x_{h}$ for $h \in {\mathcal E}$ such that 
$$
x_{e} \neq x_{-e} \ (e, -e \in {\mathcal E}), \ 
x_{h} \neq x_{h'} \ (\mbox{$h, h' \in {\mathcal E}$ with $h \neq h'$ and $v_{h} = v_{h'}$}).
$$
Then we obtain a degenerate complex curve $R_{\Delta}$ of genus $g$  
with $n$ marked points $x_{t} \in P_{v_{t}}$ $(t \in T)$ 
by identifying $x_{h} \in P_{v_{h}}$ with $x_{-h} \in P_{v_{-h}}$ $(h \in \pm E)$, 
where $x_{h} = \infty$ if $h \in {\mathcal E}_{\infty}$. 
Denote by $\xi_{h}$ $(h \in \pm E)$ denotes the local coordinate 
on $P_{v(h)} = {\mathbb P}^{1}$ given by 
$$
\xi_{h}(z) = \left\{ \begin{array}{ll} 
z - x_{h} & \mbox{(if $h \not\in {\mathcal E}_{\infty}$),} 
\\
1/z   & \mbox{(if $h \in {\mathcal E}_{\infty}$).} 
\end{array} \right. 
$$
Then for non-zero complex deformation parameters $y_{e} = y_{-e}$ $(e \in E)$, 
a universal deformation ${\mathcal R}_{\Delta}$ (with marked points $x_{t}$ $(t \in T)$) 
of $R_{\Delta}$ is defined by gluing $\bigcup_{v \in V} P_{v}$ minus neighborhoods of 
$x_{h}$ $(h \in E)$ via the relations $\xi_{h} \cdot \xi_{-h} = y_{h}$. 

We describe ${\mathcal R}_{\Delta}$ as a family of Schottky uniformized Riemann surfaces. 
Denote by $\phi_{h}$ an element of $PGL_{2}({\mathbb C})$ defined as  
$$
\xi_{h} \left( \phi_{h}(z) \right) \xi_{-h}(z) = y_{h} \ \ \left( z \in {\mathbb P}^{1} \right). 
\eqno({\rm A.1}) 
$$
We note that the attractive fixed point $\alpha_{\phi_{h}}$, 
the repulsive fixed point $\alpha'_{\phi_{h}}$ and the multiplier $\beta_{\phi_{h}}$ 
of $\phi_{h}$ are expressed by $x_{\pm h}$, $y_{h}$. 
If $h \in {\mathcal E}_{\infty}$ (resp. $-h \in {\mathcal E}_{\infty}$), 
then $\alpha_{\phi_{h}}$ is $\infty$ (resp. $x_{h}$), 
$\alpha'_{\phi_{h}}$ is $x_{-h}$ (resp. $\infty$) and $\beta_{\phi_{h}}$ is $y_{h}$. 
If $\pm h \not\in {\mathcal E}_{\infty}$, 
then using Hensel's lemma we can compute the solutions to 
the characteristic polynomial of $\phi_{h}$ 
as $x = u \cdot {\rm tr}(\phi_{h})$, $x' = \nu u^{-1} \cdot {\rm tr}(\phi_{h})$, 
where $\nu = \det(\phi_{h})/{\rm tr}(\phi_{h})^{2}$, 
and hence $\beta_{\phi_{h}} = \nu u^{-2}$ and 
$\xi_{h} \left( \alpha_{\phi_{h}} \right)$, $\xi_{-h} \left( \alpha'_{\phi_{h}} \right)$ 
are computable. 
For any reduced path $\rho = h(1) \cdot h(2) \cdots h(l)$ in $\Delta$ 
which is the product of oriented edges $h(1), ... ,h(l)$ such that $v_{h(i)} = v_{-h(i+1)}$, 
one can associate an element $\rho^{*}$ of $PGL_{2}({\mathbb C})$ 
with reduced expression $\phi_{h(l)} \phi_{h(l-1)} \cdots \phi_{h(1)}$. 
Fix a base vertex $v_{b}$ of $V$, 
and consider the fundamental group $\pi_{1} (\Delta; v_{b})$ 
which is a free group of rank $g = \dim H_{1}(\Delta)$. 
Then the correspondence $\rho \mapsto \rho^{*}$ 
gives an injective map $\pi_{1} (\Delta; v_{b}) \rightarrow PGL_{2}({\mathbb C})$ 
whose image is denoted by $\Gamma_{\Delta}$. 
Under the assumption that $|y_{e}|$ $(e \in E)$ are sufficiently small, 
$\Gamma_{\Delta}$ are Schottky groups, 
and $R_{\Gamma_{\Delta}}$ (with marked points $x_{t}$ $(t \in T)$) 
give the universal deformation ${\mathcal R}_{\Delta}$ of $R_{\Delta}$ 
(cf. \cite[Section 3]{I1}).

\subsection{Comparison of deformations} 

Assume that $\Delta = (V, E, T)$ is a stable graph which is not trivalent. 
Then there exists a vertex $v_{0} \in V$ which has at least $4$ branches. 
Take two elements $h_{1}, h_{2}$ of $\pm E \cup T$ such that 
$h_{1} \neq h_{2}$ and $v_{h_{1}} = v_{h_{2}} = v_{0}$, 
and let $\Delta' = (V', E', T')$ be a stable graph obtained from $\Delta$ 
by replacing $v_{0}$ with an oriented (non-loop) edge $h_{0}$ such that 
$v_{h_{1}} = v_{h_{2}} = v_{h_{0}}$ and that $v_{h} = v_{-h_{0}}$ 
for any $h \in \pm E \cup T - \{ h_{1}, h_{2} \}$ with $v_{h} = v_{0}$. 
Put $e_{i} = |h_{i}|$ $(i = 0, 1, 2)$. 
Then we have the following identifications: 
$$
V = V' - \{ v_{h_{0}} \} \ (\mbox{in which $v_{0} = v_{-h_{0}}$}), 
\ E = E' - \{ e_{0} \}, \ T = T'. 
$$ 
Denote by $t_{h}$ (resp. $s_{e}$) the moduli (resp. deformation) parameters of 
${\mathcal R}_{\Delta'}$ corresponding to $h \in \pm E' \cup T'$ (resp. $e \in E'$). 
We regard ${\mathcal R}_{\Delta}$ and ${\mathcal R}_{\Delta'}$ with marked points 
indexed by $T = T'$ as universal deformations of $R_{\Delta}$ 
by $y_{e}$ $(e \in E)$ and $s_{e}$ $(e \in E' - \{ e_{0} \})$ respectively. 
Denote by ${\mathbb C} [[s_{e}]]$ the ring of formal power series in $s_{e}$ $(e \in E')$, 
and by ${\mathbb C} [[s_{e}]]^{\times}$ its subset 
consisting of formal power series with non-zero constant terms. 
The ratio $\left( \phi_{-h_{0}}(t_{h_{1}}) - \phi_{-h_{0}}(t_{h_{2}}) \right)/s_{e_{0}}$ belongs to 
${\mathbb C} [[s_{e}]]^{\times}$, 
and hence by the above universality, 
${\mathcal R}_{\Delta}$ is isomorphic to ${\mathcal R}_{\Delta'}$ 
over ${\mathbb C} [[s_{e}]] \left[ e_{0}^{-1} \right]$ such that 
$$
\frac{x_{h_{1}} - x_{h_{2}}}{s_{e_{0}}}, \ \ 
\frac{y_{e_{i}}}{s_{e_{0}} s_{e_{i}}} \ (\mbox{$i = 1, 2$ with $h_{i} \not\in T$}), \ \
\frac{y_{e}}{s_{e}} \ (e \in E - \{ e_{1}, e_{2} \}) 
$$
belong to ${\mathbb C} [[s_{e}]]^{\times}$ if $h_{1} \neq - h_{2}$ 
since the multipliers of $\phi_{-h_{0}} \phi_{h_{i}}$ $(i = 1, 2)$ are equal to $s_{e_{0}} s_{e_{i}}$ 
modulo ${\mathbb C} [[s_{e}]]^{\times}$, 
and 
$$
\frac{x_{h_{1}} - x_{h_{2}}}{s_{e_{0}}}, \ \ 
\frac{y_{e}}{s_{e}} \ (e \in E) 
$$
belong to ${\mathbb C} [[s_{e}]]^{\times}$ if $h_{1} = - h_{2}$. 

By this comparison theorem, 
we obtain a family of Schottky uniformized Riemann surfaces (with marked points) 
close to degenerate complex curves. 
\vspace{2ex}

\noindent 
{\bf Acknowledgments} 
\vspace{2ex}

This work is partially supported by the JSPS Grant-in-Aid for 
Scientific Research No. 20K03516.


\begin{thebibliography}{99}

%\begin{thebibliography}{PTW02} % '2nd argument contains the widest acronym'

\bibitem{GiEG} 
J. I. B. Gill, K. Ebrahimi-Fard and H. Gangl, 
Periods in quantum field theory and arithmetic, 
ICMAT, Madrid, 2014, 
Springer Nature, 2020. pp. 630 

\bibitem{Brod}
J. Br\"{o}del, C. Duhr, F. Dulat, B. Penante and L. Tancredi, 
Elliptic symbol calculus: from elliptic polylogarithms to 
iterated integrals of Eisenstein series, 
JHEP {\bf 08} (2018) 014, arXiv:1803.10256 

\bibitem{Bo}
C. Bogner, S. M\"{u}ller-Stach and S. Weinzierl, 
The unequal mass sunrise integral expressed through iterated integrals on 
$\overline{\mathcal M}_{1,3}$, 
Nucl. Phys. B {\bf 954} (2020), 114991, arXiv:1907.01251 

\bibitem{W} 
S. Weinzierl, 
Feynman Integrals, A Comprehensive Treatment for Students and Researchers, 
UNITEXT for Physics (2022), arXiv:2201.03593 

\bibitem{Bou}
J. L. Bourjaily, J. Broedel, E. Chaubey, C. Duhr, H. Frellesvig, M. Hidding, R. Marzucca, 
A. J. McLeod, M. Spradlin, L. Tancredi, C. Vergu, M. Volk, A. Volovich, M. v. Hippel, 
S. Weinzierl, M. Wilhelm and C. Zhang, 
Functions beyond multiple polylogarithms for precision collider physics, 
arXiv:2203.07088 

\bibitem{DHS}
E. D'Hoker, M. Hidding and O. Schlotterer, 
Constructing polylogarithms on higher-genus Riemann surfaces, 
arXiv:2306.08644 (2023) 

\bibitem{CEE} 
D. Calaque, B. Enriquez and P. Etingof, 
Universal KZB equations: the elliptic case, 
in Algebra, arithmetic, and geometry: in honor of Yu. I. Manin. Vol. I, 165--266, 
Progr. Math., 269, Birkh\"{a}user, Boston, 2009, 
arXiv: math/0702670 

\bibitem{E1} 
B. Enriquez, 
Elliptic associators, 
Sel. Math. New Ser. {\bf 20} (2014), 491--584. 

\bibitem{H} 
R. Hain, 
Notes on the universal elliptic KZB equation, 
Pure and Applied Mathematics Quarterly, vol. 12, no. 2 (2016) International Press, 
arXiv: 1309.0580 

\bibitem{D2} 
P. Deligne, 
Le groupe fondamental de la droite projective moins trois points, 
in Galois groups over ${\mathbb Q}$, 
Publ. MSRI, vol. 16, Springer, 1989, pp. 79--298. 

\bibitem{Brow} 
F. Brown, 
Multiple zeta values and periods of moduli spaces $\overline{\mathfrak M}_{0,n}$, 
Ann. Sci. \'{E}c. Norm. Sup\'{e}r. {\bf 42} (2009), 371--489. 

\bibitem{BPP} 
P. Banks, E. Panzer and B. Pym, 
Multiple zeta values in deformation quantization, 
Invent. Math. https://doi.org/10.1007/s00222-020-00970-x, 
arXiv: 1812.11649 

\bibitem{GoM} 
A. B. Goncharov and Y. I. Manin, 
Multiple $\zeta$-motives and moduli spaces $\overline{\mathcal M}_{0,n}$, 
Compos. Math. {\bf 140} (2004), 1--14. 

\bibitem{E2} 
B. Enriquez, 
Analogues elliptiques des nombres multiz\'{e}tas, 
arXiv:1301.3042

\bibitem{DM} 
P. Deligne and D. Mumford, 
The irreducibility of the space of curves of given genus, 
Publ. Math. IHES {\bf 36} (1969), 75--109. 

\bibitem{K} 
F. F. Knudsen, 
The projectivity of the moduli space of stable curves II, III, 
Math. Scand. {\bf 52} (1983), 161--199, 200--212. 

\bibitem{KM} 
F. F. Knudsen and D. Mumford, 
The projectivity of the moduli space of stable curves I, 
Math. Scand. {\bf 39} (1976), 19--55. 

\bibitem{I1} 
T. Ichikawa, 
Generalized Tate curve and integral Teichm\"{u}ller modular forms, 
Amer. J. Math. {\bf 122} (2000), 1139--1174. 

\end{thebibliography}
\end{document}